\documentclass[reqno,12pt]{amsart}

\usepackage{epsf}
\usepackage{graphics}
\usepackage{amssymb}
\usepackage{amsmath}

\date{}

\theoremstyle{plain}
\newtheorem{theorem}{Theorem}

\theoremstyle{definition}

\theoremstyle{remark}

\title{The volume of positive braid links} 

\author{Sebastian Baader}

\begin{document}

\begin{abstract} Based on recent work by Futer, Kalfagianni and Purcell, we prove that the volume of sufficiently complicated positive braid links is proportional to the signature defect $\Delta \sigma=2g-\sigma$.
\end{abstract}

\maketitle

\section{Introduction}

A positive braid on $n$ strings is a product of positive standard generators of the braid group $B_n$. The canonical closure of a positive braid $\beta \in B_n$ is a link $\hat{\beta}$ with at most $n$ components.
Let $b_1$ be the first Betti number of the canonical fibre surface of the link $\hat{\beta}$. Furthermore, let $\sigma$ be the signature invariant of the link $\hat{\beta}$. The quantity $\Delta \sigma(\hat{\beta})=b_1-\sigma$ is called the signature defect of the link $\hat{\beta}$. In the case of knots, this is simply $2g-\sigma$, where $g$ is the minimal genus of $\hat{\beta}$. The main result of this note is a volume estimate in terms of the signature defect.

\begin{theorem} Let $\hat{\beta} \subset S^3$ be a hyperbolic link associated with a sufficiently complicated positive braid $\beta$. Then
$$\frac{1}{3} v_8 \, \Delta \sigma(\hat{\beta}) \leq \text{vol}(S^3 \setminus \hat{\beta}) < 105v_3 \, \Delta \sigma(\hat{\beta}),$$
where $v_3=1.0149...$ and $v_8=3.6638...$ are the volumes of a regular ideal tetrahedron and octahedron, respectively.
\end{theorem}

Here a positive braid is sufficiently complicated, if it is a product of powers $\sigma_i^k$, with $k \geq 3$. The twist number $t(\beta)$ of a sufficiently complicated braid $\beta$ is the minimal number of factors in such a product.
As discussed in~\cite{FKP}, the link $\hat{\beta}$ is hyperbolic, if and only if $\beta$ is prime. Equivalently, $\beta$ contains at least two non-consecutive factors of the form $\sigma_i^k$, for all $i \leq n-1$. In the recent monograph~\cite{FKP}, Futer, Kalfagianni and Purcell determined tight bounds for the volume of various families of hyperbolic links, in terms of the twist number. In particular, they proved the following volume estimates for sufficiently complicated positive braid links.

\begin{theorem}[\cite{FKP}, Theorem~9.7]
Let $\hat{\beta} \subset S^3$ be a hyperbolic link associated with a sufficiently complicated positive braid $\beta$. Then
$$\frac{2}{3} v_8 \, t(\beta) \leq \text{vol}(S^3 \setminus \hat{\beta}) < 10 v_3 \, (t(\beta)-1).$$
\end{theorem}

\noindent
Theorem~1 is an immediate consequence of this and the following estimates for the twist number, which we will prove in the next section.

\begin{theorem} Let $\hat{\beta} \subset S^3$ be a hyperbolic link associated with a sufficiently complicated positive braid $\beta$. Then
$$\frac{1}{2} \, \Delta \sigma(\hat{\beta}) \leq t(\beta) \leq \frac{21}{2} \, \Delta \sigma(\hat{\beta}).$$
\end{theorem}

\section*{Acknowledgements}

I would like to thank Julien March\'e for drawing my attention to the work of Futer, Kalfagianni and Purcell.

\section{Twist number and signature defect}

The signature $\sigma(L)$ of a link $L$ is defined as the signature of any symmetrised Seifert matrix $V+V^T$ associated with $L$. Let $b_1$ be the first Betti number of a minimal genus Seifert surface for $L$. Then
$$-b_1 \leq \sigma(L) \leq b_1.$$
In particular, the signature defect $\Delta \sigma=b_1-\sigma$ is positive. The proof of Theorem~3 relies on three elementary facts about the signature invariant.

\medskip
\begin{enumerate}
\item The signature defect of the closure of $\sigma_1^n$ is zero.\\
\item The signature defect of the closure of $\sigma_1^{k_1} \sigma_2^{k_2} \sigma_1^{k_3} \sigma_2^{k_4}$ is two, provided all $k_i \geq 2$ (this fact is a central ingredient in the classification of positive braids with $\Delta \sigma=0$, see~\cite{Ba}).\\
\item Let $\Sigma \subset \widetilde{\Sigma}$ be an inclusion of Seifert surfaces which induces an injection on the level of first homology groups. Then
$$\Delta \sigma(\partial \widetilde{\Sigma}) \leq \Delta \sigma(\partial \Sigma)+2(b_1(\widetilde{\Sigma})-b_1(\Sigma)).$$
\end{enumerate}

\begin{proof}[Proof of Theorem~3]
Let $\beta$ be a sufficiently complicated positive braid with twist number $t(\beta)$; let $\widetilde{\Sigma} \subset S^3$ be the canonical fibre surface of the link $\hat{\beta}$ (see Stallings~\cite{St}). By cutting the surface $\widetilde{\Sigma}$ along an interval on the left of every twist region $\sigma_i^k$, as sketched in Figure~1, we obtain a subsurface $\Sigma \subset \widetilde{\Sigma}$.

\begin{figure}[ht]
\scalebox{1.0}{\raisebox{-0pt}{$\vcenter{\hbox{\epsffile{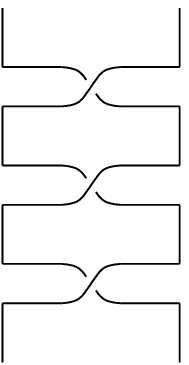}}}$}} \qquad $\longrightarrow$ \qquad
\scalebox{1.0}{\raisebox{-0pt}{$\vcenter{\hbox{\epsffile{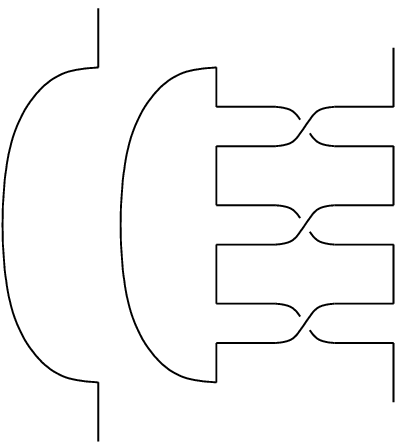}}}$}}
\caption{}
\end{figure}

The boundary of $\Sigma$ is a disjoint union of connected sums of torus links on two strings. According to fact~(1), the signature defect $\Delta \sigma(\partial \Sigma)$ is zero. Using fact~(3), we conclude
$$\Delta \sigma(\hat{\beta})=\Delta \sigma(\partial \widetilde{\Sigma}) \leq 2(b_1(\widetilde{\Sigma})-b_1(\Sigma)) \leq 2t(\beta).$$
This is the first inequality of Theorem~3.

\medskip
For the second inequality, we first consider the case of 3-string braids. Let $\beta \in B_3$ be a sufficiently complicated positive braid with hyperbolic closure. Then $t(\beta) \geq 4$. Moreover, $\beta$ contains at least $\frac{t(\beta)}{7}$ consecutive subwords of the form
$$\sigma_1^{k_1} \sigma_2^{k_2} \sigma_1^{k_3} \sigma_2^{k_4}$$
(a better estimate would be $\frac{1}{4}(t(\beta)-3)$, but this is not linear). By fact~(2), every such subword contributes two to the signature defect. More precisely, the fibre surface of $\hat{\beta}$ contains at least $\frac{t(\beta)}{7}$ disjoint subsurfaces whose homology groups are orthogonal with respect to the symmetrised Seifert form. All these contribute two to the signature defect. In particular,
$$\frac{2}{7} t(\beta) \leq \Delta \sigma(\hat{\beta}).$$

Now let us turn to the case of higher braid indices. Our goal is to find a large number of non-interfering subwords of the form $\sigma_i^{k_1} \sigma_{i+1}^{k_2} \sigma_i^{k_3} \sigma_{i+1}^{k_4}$ in $\beta$. For this purpose, let us partition the strings of the braid $\beta \in B_n$ into three subsets $S_1$, $S_2$, $S_3$, according to their index modulo $3$. A simple counting argument shows that one of the subsets $S_j$ carries at least $\frac{t(\beta)}{21}$ disjoint subwords of the desired form. Here `carrying' means that the central string $i+1$ of the word $\sigma_i^{k_1} \sigma_{i+1}^{k_2} \sigma_i^{k_3} \sigma_{i+1}^{k_4}$ belongs to $S_j$. Indeed, let us put dots on the strings of the closed braid $\hat{\beta}$, one between each pair of adjacent twist regions, as in $\sigma_i^k \sigma_{i+1}^l$ or $\sigma_i^k \sigma_{i-1}^l$. Altogether, there are at least $t(\beta)$ dots, otherwise the number of factors $\sigma_i^k$ would not be minimal. One of the subsets $S_j$ carries at least $\frac{t(\beta)}{3}$ dots, thus at least $\frac{t(\beta)}{3 \cdot 7}$ disjoint subwords of the form $\sigma_i^{k_1} \sigma_{i+1}^{k_2} \sigma_i^{k_3} \sigma_{i+1}^{k_4}$. Here it is important that neighbouring strings of $S_j$ are at distance three. As before, we conclude
$$\frac{2}{21} t(\beta) \leq \Delta \sigma(\hat{\beta}).$$

\end{proof}

\bigskip
\noindent
Universit\"at Bern, Sidlerstrasse 5, CH-3012 Bern, Switzerland

\bigskip
\noindent
\texttt{sebastian.baader@math.unibe.ch}


\begin{thebibliography}{9}

\bibitem{Ba}
     S.~Baader: \emph{Positive braids of maximal signature}, arXiv:1211.4824. 

\bibitem{FKP}
     D.~Futer, E.~Kalfagianni, J.~Purcell: \emph{Guts of surfaces and the colored Jones polynomial}, Lecture Notes in Mathematics, 2069. Springer, Heidelberg, 2013.
     
\bibitem{St}
     J.~Stallings: \emph{Constructions of fibred knots and links}, Algebraic and geometric topology, Proc. Sympos. Pure Math.~\textbf{32} (1978), 55-60, Amer. Math. Soc., Providence, R.I.

\end{thebibliography}
\end{document}